\documentclass[12pt,draft]{extartic}

\usepackage{amsmath,amsthm,amssymb,calc}
\usepackage[dvips]{graphics, graphicx}

\usepackage[english]{babel}

\newcounter{num}[section]

\newcommand{\Th}{\refstepcounter{num}
{\bf Theorem \arabic{section}.\arabic{num} }}

\newcommand{\Lemma}{\refstepcounter{num}
{\bf Lemma \arabic{section}.\arabic{num} }}

\newcommand{\Pred}{\refstepcounter{num}
{\bf Proposition \arabic{section}.\arabic{num} }}

\newcommand{\Cor}{\refstepcounter{num}
{\bf Corollary \arabic{section}.\arabic{num} }}

\newcommand{\Note}{\refstepcounter{num}
{\it Note \arabic{section}.\arabic{num} }}

\newcommand{\St}{\refstepcounter{num}
{\bf Statement \arabic{section}.\arabic{num} }}

\newcommand{\Def}{\refstepcounter{num}
{\it Definition \arabic{section}.\arabic{num} }}

\newcommand{\Proof}{{\bf Proof. }}

\def\_phi{\varphi}

\def\a{\alpha}

\def\la{\lambda}

\def\v{\vec}
\def\F{\widehat}

\def\m{\times}

\def\C{{\mathbb C}}

\def\Z_N{{\mathbb Z}_N}
\def\Z{{\mathbb Z}}

\def\f{{\mathbb F}}
\def\Gr{{\mathbf G}}

\def\l{\left}
\def\r{\right}

\def\supp{{\rm supp\,}}

\author{Shkredov I.D.\footnote{
The author is supported
Pierre Deligne's grant based on his 2004 Balzan prize, President's
of Russian Federation grant N МК--1959.2009.1, grant RFFI N
06-01-00383 and grant Leading Scientific Schools No. 8684.2010.1},
Vyugin I.V.\footnote{The author is supported by grant Leading Scientific Schools No. 8508.2010.1}}
\title{On additive shifts of multiplicative subgroups}
\date{}
\begin{document}
\maketitle

\begin{center}
    Annotation.
\end{center}

{\it \small

    Generalizing a result of S.V. Konyagin and D.R. Heath--Brown, we prove, in particular,
    that for any multiplicative subgroup $R \subseteq \Z/p\Z$
    and
    any nonzero elements
    $\mu_1,\dots,\mu_k$
    the following holds
    $|R\bigcap (R+\mu_1) \bigcap \dots \bigcap (R+\mu_k)| \ll_k |R|^{\frac{1}{2}+\alpha_k}$,
    provided by $1 \ll_k |R| \ll_k p^{1-\beta_k}$,
    where $\alpha_k$, $\beta_k$ are some sequences of positive reals and $\alpha_k, \beta_k \to 0$, $k\to \infty$.
    Besides we show that for an arbitrary subgroup $R$, $|R| \ll p^{1/2}$
    one have $|R\pm R| \gg |R|^{5/3} \log^{-1/2} |R|$.
}
\\

\refstepcounter{section}
\label{sec:introduction}

{\bf \arabic{section}. Introduction.}

Let $p$ be a prime number, $\Z_p^* = (\Z/p\Z) \setminus \{ 0 \}$
be the group of all  invertible elements of the field $\Z_p$, and
$R\subseteq \Z_p^*$ be its multiplicative subgroup.
Different properties of such subgroups have been studied by several authors, see e.g.
\cite{Bou_prod1}--\cite{BK},
\cite{Waring_Z_p}--\cite{Schoen_Shkr},
\cite{Shkredov_RplusR},
\cite{Yekhanin_subgroups}.
For example A. Garcia and J.F. Voloch \cite{Garcia_Voloch}, using deep algebraic ideas,
proved that for any subgroup $R$, $|R|< (p-1)/ ((p-1)^{1/4}+1)$
and an arbitrary nonzero $\mu$
the following holds
\begin{equation}\label{f:Garcia_Voloch}
    |R \bigcap (R+\mu)| \le 4 |R|^{2/3} \,.
\end{equation}
D.R. Heath--Brown and  S.V. Konyagin generalized (\ref{f:Garcia_Voloch}) and gave another prove of the result
in \cite{Heath_B-K}
(see also \cite{K_Tula}).
Their approach uses a well--known method of S.A. Stepanov \cite{Stepanov}.
In the paper we extend the result of Garcia--Voloch and also similar theorems from
\cite{Heath_B-K}, \cite{K_Tula} for the case of several additive shifts.
Let us formulate one of the main of our results.

\Th
{\it
    Let $R\subseteq \Z_p^*$ be a multiplicative subgroup, $k\ge 1$ be a positive integer,
    $|R| > k 2^{2k+4}$.
    Let also $\mu_1,\dots,\mu_k$ be different nonzero residuals, and $Q=RQ$ be a $R$---invariant set,
    $0\notin Q$,
    $|Q| < ( (|R|/k)^{1/2k} - 1)^{2k+1}$, $p\ge 4k|R| ( |Q|^{\frac{1}{2k+1}} + 1)$.
    Then
    \begin{equation}\label{f:main_t}
        \sum_{\la \in Q} |R \bigcap (R +\la \cdot \mu_1) \bigcap \dots \bigcap (R +\la \cdot \mu_k)|
            \le
                4 (k+1) (|Q|^{\frac{1}{2k+1}} + 1)^{k+1} |R| \,.
    \end{equation}
}
\label{t:main_many_shifts}

Theorem \ref{t:main_many_shifts} easily implies a statement on
the maximal cardinality of the intersection of
$k$
additive shifts of a subgroup.

\Cor
{\it
    Let $R\subseteq \Z_p^*$ be a multiplicative subgroup,
    $k\ge 1$ be a positive integer, and $\mu_1,\dots,\mu_k$ be different nonzero elements.
    Let also
    $$
        32 k 2^{20k \log (k+1)} \le |R|\,, \quad  p \ge 4k |R|  ( |R|^{\frac{1}{2k+1}} + 1 ) \,.
    $$
    Then
    $$
        |R\bigcap (R+\mu_1) \bigcap \dots (R+\mu_k)| \le 4 (k+1) (|R|^{\frac{1}{2k+1}} + 1)^{k+1} \,.
    $$
}
\label{cor:main_many_shifts}

Roughly speaking, the corollary above asserts that
$|R\bigcap (R+\mu_1) \bigcap \dots \bigcap (R+\mu_k)| \ll_k |R|^{\frac{1}{2}+\alpha_k}$,
provided by
$1 \ll_k |R| \ll_k p^{1-\beta_k}$,
where $\alpha_k, \beta_k$ are some sequences of positive numbers, and $\alpha_k, \beta_k \to 0$, $k\to \infty$.

Our approach develops the method from \cite{Heath_B-K}, \cite{K_Tula}.

Now consider another additive characteristic of multiplicative subgroups,
namely, the cardinality of their sums and differences.
Bound (\ref{f:Garcia_Voloch}) implies that (see \cite{Garcia_Voloch})
$$
    |R\pm R| \gg |R|^{4/3}
$$
for any subgroup $R$ with $|R| \ll p^{3/4}$.
D.R. Heath--Brown and  S.V. Konyagin in \cite{Heath_B-K} (see also \cite{K_Tula})
proved
\begin{equation}\label{f:HK_3/2}
    |R\pm R| \gg |R|^{3/2}
\end{equation}
for all subgroups $R$ such that $|R| \ll p^{2/3}$.
Using a combinatorial idea from \cite{Katz_Koester}
(see also papers \cite{Sanders_non-abealian}---\cite{Schoen_Freiman}, which are develop the approach),
we improve inequality (\ref{f:HK_3/2})
(see Theorem \ref{t:subgroups_doubling} of section \ref{sec:applications})
in the following way
\begin{equation}\label{f:5/3_introduction}
    |R\pm R| \gg \frac{|R|^{5/3}}{\log^{1/2} |R|}
\end{equation}
for subgroups $R$ with the condition $|R| \ll p^{1/2}$.

Let us say a few words about the structure of the paper.
In auxiliary section \ref{sec:Katz-Koester} we give a series of required definitions
and discuss, in detail, a generalization of ordinary convolutions,
which is naturally appears in the problems concerning several additive shifts.
In the next section \ref{sec:Vronskian} we obtain preliminary results
on linear dependence of some systems of polynomials in $\Z_p [x]$.
Applying Stepanov's method and using linear independence of such polynomials,
we get
Theorem \ref{t:main_many_shifts}
in the next section \ref{sec:proof}.
The last section \ref{sec:applications} contains consequences
of the obtained results,
and also their applications to combinatorial number theory.
Here we prove, in particular,
inequality
(\ref{f:5/3_introduction}).

We conclude with few comments regarding the notation used in this paper.
Let $\Z_p = \Z/p\Z$, and $\Z_p^* = \Z_p \setminus \{ 0 \}$.
If $A$ is a set then we write $A(x)$ for its characteristic function.
Thus $A(x) = 1$ if $x\in A$ and $A(x)=0$ otherwise.
We use the symbol $|A|$ to denote the cardinality of the set $A$.
All logarithms $\log$ are base $2.$
Signs $\ll$ and $\gg$ are the usual Vinogradov's symbols.
For a positive integer $n,$ we set $[n]=\{1,\ldots,n\}.$

The authors are grateful of S.V. Konyagin for a number of helpful advices and remarks.

\refstepcounter{section}
\label{sec:Katz-Koester}

{\bf \arabic{section}. Katz--Koester method and higher convolutions.}

Recall the required definitions.
Let $\Gr$ be a finite Abelian group, $N=|\Gr|.$
It is well--known~\cite{Rudin_book} that the dual group $\F{\Gr}$ is isomorphic to $\Gr.$
Let $f$ be a function from $\Gr$ to $\C.$  We denote the Fourier transform of $f$ by~$\F{f},$
\begin{equation}\label{F:Fourier}
  \F{f}(\xi) =  \sum_{x \in \Gr} f(x) e( -\xi \cdot x) \,,
\end{equation}
where $e(x) = e^{2\pi i x}$.
Define the two convolutions of functions $f$ and $g$
$$
    (f*g) (x) := \sum_{y\in \Gr} f(y) g(x-y)\,, \quad \mbox{ and } \quad (f\circ g) (x) := \sum_{y\in \Gr} f(y) g(y+x) \,.
$$
Write $E(A,B)$ for {\it additive energy} of two sets $A,B \subseteq \Gr$ (see e.g. \cite{Tao_Vu_book}), that is
$$
    E(A,B) = |\{ a_1+b_1 = a_2+b_2 ~:~ a_1,a_2 \in A,\, b_1,b_2 \in B \}| \,.
$$
If $A=B$ we simply write $E(A)$ instead of $E(A,A).$ Clearly,
\begin{equation}\label{f:energy_convolution}
    E(A,B) = \sum_x (A*B) (x)^2 = \sum_x (A \circ B) (x)^2 = \sum_x (A \circ A) (x) (B \circ B) (x)
    \,.
\end{equation}

Consider a generalization of the operation $\circ$.

\Def
\label{def:convolution}
{
    Let $k\ge 1$ be a positive number, and $f_1,\dots,f_k : \Gr \to \C$ be functions.
Denote by $C_k (f_1,\dots,f_k) (x_1,\dots, x_{k-1})$ the function
$$
    C_k (f_1,\dots,f_k) (x_1,\dots, x_{k-1}) = \sum_z f_1 (z) f_2 (z+x_1) \dots f_k (z+x_k) \,.
$$
Thus, $C_2 (f_1,f_2) (x) = (f_1 \circ f_2) (x)$.
Put $C_1 (f) = \sum_z f(z)$.
If $f_1=\dots=f_k=A$, $A\subseteq \Gr$ is a set then write
$C_k (A) (x_1,\dots, x_{k-1})$ for $C_k (f_1,\dots,f_k) (x_1,\dots, x_{k-1})$.
}

\Def
\label{def:tensor_product}
{
    Let $A,B \subseteq \Gr$ be arbitrary sets and $l\ge 1$ be a positive integer.
    Then
    \begin{equation}\label{def:otimes}
        A \otimes_l B = \bigcup_{b\in B} (A-b)^l \subseteq \Gr^l \,.
    \end{equation}
    In particular $A \otimes_1 B = A \otimes B = A-B$.
}

Clearly,
$$
    \supp C_k (B,A,\dots,A) = \bigcup_{a\in B} (A-a)^{k-1} = A \otimes_{k-1} B \subseteq \Gr^{k-1} \,.
$$
We have  $|A|^{k-1} \le |A \otimes_{k-1} B| \le |B| |A|^{k-1}$.
In particular, the set $A \otimes_{k-1} B$ is nonempty.
Let
$$
    E_k (f_1,\dots,f_k) = \sum_{x_1,\dots,x_{k-1}} C^2_k (f_1,\dots,f_k) (x_1,\dots, x_{k-1}) \,.
$$
Then $E_2(A,B) = E(A,B)$.
We write $E_k (A)$ for $E_k (A,\dots,A)$.
There is an obvious connection between quantities  $|A \otimes_{k-1} A|$ and $E_k (A)$.

\Lemma
{\it
    Let $A,B\subseteq \Gr$ be two sets, and $k\ge 2$ be a positive integer.
    Then
    $$
        |A|^{2k-2} |B|^2 \le E_k (A,\dots,A,B) \cdot |A \otimes_{k-1} B| \,.
    $$
}
\label{l:KB_high_convolutions}
\Proof
We have $\sum_{x_1,\dots,x_{k-1}} C_k (B,A,\dots,A) (x_1,\dots, x_{k-1}) = |A|^{k-1} |B|$.
Using Cauchy--Schwarz, we obtain the required estimate.
$\Box$

Let $B\subseteq A$ be a set, and $(x_1,\dots,x_k) := \v{x} \in A \otimes_{k-1} B$
be a vector.
Put $B_{\v{x}} = B \bigcap (A-x_1) \bigcap (A-x_2) \bigcap \dots \bigcap (A-x_k)$.
Clearly, $B_{\v{x}}$ is nonempty.
Besides $|B_{\v{x}}| = C_k (B,A,\dots,A) (x_1,\dots,x_{k-1})$.
We can easily describe the structure $A \otimes_{k-1} B$ using the sets  $B_{\v{x}}$.

\Lemma
{\it
    Let $B\subseteq A \subseteq  \Gr$ be two sets,
    and $l\ge 1$ be a positive integer.
    Then
    $$
        A \otimes_{l} B = \{ (x_1,\dots,x_{l}) ~:~ A_{x_1,\dots,x_{l}} \bigcap B \neq \emptyset \} \,.
    $$
}

\Cor
{\it
    Let $B\subseteq A \subseteq  \Gr$ be two sets,
    and $l\ge 2$, $m\ge 1$ be positive integers, $m\le l$.
    Then
    \begin{equation}\label{f:A_l-A_m_eq_1}
        A \otimes_l B = \bigcup_{(x_1,\dots,x_m) \in A \otimes_m B} \, \{ (x_1,\dots,x_m) \} \times (A \otimes_{l-m} B_{x_1,\dots,x_m}) \,.
    \end{equation}
    In particular,
    \begin{equation}\label{f:A_l-A_m_eq_2}
        A \otimes_l A
            =
                \bigcup_{(x_1,\dots,x_{l-1}) \in A \otimes_{l-1} A} \, \{ (x_1,\dots,x_{l-1}) \}
                    \m
                        (A-A_{x_1,\dots,x_{l-1}}) \,.
    \end{equation}
}
\label{cor:A_l-A_m_eq}

We need in upper bounds for the cardinality of $A \otimes_{k-1} A$.
For positive integers $l$ and $m$, $m\le l$,
arbitrary set $E\subseteq [l]$, $E=\{ j_1,\dots, j_m \}$,
and any vector $x=(x_1,\dots,x_l)$ the symbol $x^E$ denotes the vector $(x_{j_1}, \dots, x_{j_m})$.
The following lemma is a consequence of the definitions.

\Lemma
{\it
    Let $A \subseteq  \Gr$ be a set,
    and $l\ge 1$ be a positive integer.
    Let also $S=A-A$.
    Then
    \begin{equation}\label{f:A_l-A_m}
        (A \otimes_l A) (x_1, \dots, x_l) \le \prod_{E\subseteq [l],\, |E|=m} (A\otimes_m A) (x^E) \,.
    \end{equation}
    Besides
    \begin{equation}\label{f:A_2-S}
        (A \otimes_2 A) (x,y) \le S(x) S(y) S(x-y) \,,
    \end{equation}
    and
    \begin{equation}\label{f:A_l-S}
        (A \otimes_l A) (x_1, \dots, x_l) \le \prod_{i,j=0,\, i\neq j}^k S(x_i-x_j) \,,
    \end{equation}
    where $x_0$ denotes $0$.
}
\label{l:pre_card_KK}

Clearly,  (\ref{f:A_l-S})
is a consequence of  (\ref{f:A_l-A_m}) and (\ref{f:A_2-S}).

\Cor
{\it
    Let  $A \subseteq  \Gr$ be a set,
    and $l\ge 1$ be a positive integer.
    Let also $S=A-A$.
    Then
    \begin{equation}\label{f:card_KK}
        |A\otimes_l A| \le \sum_{x\in S} |A-A_x|^{l-1} \le \sum_{x\in S} (S\circ S)^{l-1} (x) \,.
    \end{equation}
}
\label{cor:card_KK}
\Proof
The first inequality in  (\ref{f:card_KK}) is a consequence of formula (\ref{f:A_l-A_m_eq_1})
of Corollary \ref{cor:A_l-A_m_eq}, applying with $m=1$.
Lemma \ref{l:pre_card_KK} immediately implies
the bound $|A\otimes_l A| \le \sum_{x\in S} (S\circ S)^{l-1} (x)$.
Finally, the middle inequality is a consequence of Katz--Koester inclusion \cite{Katz_Koester}
\begin{equation}\label{f:Katz-Koester_gen}
    A_{\v{s}} - A_{\v{t}} \subseteq S_{\v{u}} \,,
\end{equation}
where $\v{s} = (s_1,\dots,s_m)$, $\v{t} = (t_1,\dots,t_n)$
are two arbitrary vectors of the lengths $m,n$, respectively,
$\v{s} \in A\otimes_m A$, $\v{t} \in A\otimes_n A$, and the vector $\v{u}$
has the length $(n+1)(m+1)-1$
and consists of all non--zero sums $s_i+t_j$, $i= 0,1,\dots, m$, $j= 0,1,\dots, n$.
$\Box$.


Let us generalize Lemma 3.1 from \cite{Schoen_Shkr}.

\Lemma
{\it
    Let $A\subseteq \Gr$ be a set, $l\ge 1$, $k\ge 2$ be positive integers.
    Then
    \begin{equation}\label{f:moments_of_gen_convolutions}
        \sum_{\v{s}_1,\dots, \v{s}_k} \sum_{z_1,\dots,z_{k-1}} C^l_k (A_{\v{s}_1}, \dots, A_{\v{s}_k}) (z_1,\dots,z_{k-1})
            =
                \sum_{x_1,\dots,x_{l-1}} C^{\| \v{s} \|+k}_l (A) (x_1,\dots,x_{l-1}) \,,
    \end{equation}
    where $\| \v{s} \| = \sum_{j=1}^k |\v{s}_j|$.
    In particular,
    \begin{equation}\label{f:gen_conv1}
        \sum_{\v{s}_1,\dots, \v{s}_k} \sum_{z_1,\dots,z_{k-1}} C_k (A_{\v{s}_1}, \dots, A_{\v{s}_k}) (z_1,\dots,z_{k-1})
            = |A|^{\| \v{s} \| + k} \,,
    \end{equation}
    and
    \begin{equation}\label{f:gen_conv2}
        \sum_{\v{s}_1,\dots, \v{s}_k} \sum_{z_1,\dots,z_{k-1}} C^2_k (A_{\v{s}_1}, \dots, A_{\v{s}_k}) (z_1,\dots,z_{k-1})
            =
                \sum_{\v{s}_1,\dots, \v{s}_k} E_k (A_{\v{s}_1}, \dots, A_{\v{s}_k}) = E_{\| \v{s} \| + k} (A) \,.
    \end{equation}
}
\label{l:moments_of_gen_convolutions}
\Proof
We have (recall that $z_0=0$)
\begin{equation}\label{f:tmp_17_12_2010_1}
    \sum_{\v{s}_1,\dots, \v{s}_k} \sum_{z_1,\dots,z_{k-1}} C^l_k (A_{\v{s}_1}, \dots, A_{\v{s}_k}) (z_1,\dots,z_{k-1})
        =
            \sum_{\v{s}_1,\dots, \v{s}_k} \sum_{z_1,\dots,z_{k-1}}
                \sum_{w_1,\dots,w_l} \prod_{j=1}^l \prod_{i=1}^k A_{\v{s}_i} (w_j+z_{i-1})
\end{equation}
$$
    =
        \sum_{w_1,\dots,w_l} \sum_{z_1,\dots,z_{k-1}} C^{\| \v{s} \|}_k (A) (w_2-w_1,\dots,w_l-w_1)
            \prod_{j=1}^l \prod_{i=1}^k A (w_j+z_{i-1})
            =
$$
$$
            =
            \sum_{w_1,\dots,w_l} C^{\| \v{s} \|+k-1}_k (A) (w_2-w_1,\dots,w_l-w_1) A(w_1) \dots A(w_l)
                =
        \sum_{x_1,\dots,x_{l-1}} C^{\| \v{s} \|+k}_l (A) (x_1,\dots,x_{l-1})
        \,,
$$
because each component of any vector $\v{s}_i$ appears at formula (\ref{f:tmp_17_12_2010_1}) exactly $l$ times.
This completes the proof.  $\Box$


\refstepcounter{section}
\label{sec:Vronskian}

{\bf \arabic{section}. On linear independence of a system of polynomials.}

In paper \cite{K_Tula} the following lemma was proved.

\Lemma
{\it
    Let $\a_1 \in \Z_p^*$ be an arbitrary residual.
    Let also $t$, $B$, $D$
    be some positive integers, $p$ be a prime number, and
    \begin{equation}\label{f:polyn_parameters_pre}
        t \ge BD \,, \quad
        p \ge tB \,.
    \end{equation}
Then the polynomials of the form
\begin{eqnarray}\label{polyn_pre}
    x^{a_i} x^{tb_{0,i}}(x-\alpha_1)^{tb_{1,i}}
\end{eqnarray}
where $a_i < D$, $b_{0,i},b_{1,i} < B$ are linearly independent over $\Z_p$.
}
\label{l:sparse_pol}

In the section we generalize the lemma above for systems of polynomials with larger number of monomials.
Our dependence between parameters worse than in Lemma \ref{l:sparse_pol}.

We use the notion of formal derivative in $\mathbb{Z}_p$.
The derivative of a polynomial is
a formal derivative of the sum of its monomials, that is another polynomial
$$
\left(\sum_{i=0}^n c_ix^i\right)'=\sum_{i=1}^n ic_ix^{i-1}.
$$
We consider the derivatives of polynomials with the degree at most $p-1$.
Leibniz's law holds for the formal derivative of such polynomials.
Note that the derivation is well--defined for formal sums not functions.

\Pred
{\it
Let $n$, $t$, $B$, $D$
be positive integers, and  $p$ be a prime number.
Let also $\a_1,\dots,\alpha_n \in \Z^*_p$ be different nonzero residuals, and
\begin{equation}\label{f:polyn_parameters}
    t \ge \frac{1}{2}(n-1) B^{2n} + D B^n \,, \quad
    p \ge (2nB+2) t\,.
\end{equation}
Then the polynomials of the form
\begin{eqnarray}\label{polyn}
x^{a_i} x^{tb_{0,i}}(x-\alpha_1)^{tb_{1,i}}\ldots (x-\alpha_n)^{tb_{n,i}}
\end{eqnarray}
where $a_i < D$, $b_{0,i},b_{1,i},\ldots,b_{n,i}<B$ are linearly independent over $\Z_p$.
}
\label{pred:polyn}

\Proof
Suppose that there is a nontrivial linear combination of the polynomials from (\ref{polyn}),
which equals zero identically
\begin{eqnarray}\label{lin-komb}
\sum_{i=1}^m C_i x^{a_i} x^{tb_{0,i}}(x-\alpha_1)^{tb_{1,i}}\ldots (x-\alpha_n)^{tb_{n,i}}\equiv 0.
\end{eqnarray}
Divide  (\ref{lin-komb}) by $(x-\alpha_n)^{ts}$, where $s=\min_{i}b_{n,i}$.
Consider the terms from (\ref{lin-komb})
with minimal $b_{n,i}$, i.e. equal $s$.
One can suppose that these are the first $l_0$ terms.
Then the polynomial
\begin{eqnarray}\label{lin-komb2}
    \Phi(x)=\sum_{i=1}^{l_0} C_i x^{a_i} x^{tb_{0,i}}(x-\alpha_1)^{tb_{1,i}}\ldots (x-\alpha_{n-1})^{tb_{n-1,i}}
\end{eqnarray}
divided by $(x-\alpha_n)^t$.
Denote the sum of polynomials from
(\ref{lin-komb2}) with the same multiplier
$x^{tb_{0,i}}(x-\alpha_1)^{tb_{1,i}}\ldots (x-\alpha_{n-1})^{tb_{n-1,i}}$
as
\begin{eqnarray*}
\Phi_i(x)= H_i(x) x^{tb_{0,i}}(x-\alpha_1)^{tb_{1,i}}\ldots (x-\alpha_{n-1})^{tb_{n-1,i}},\qquad i=1,\ldots,l.
\end{eqnarray*}
\bigskip
Clearly,  $\deg H_i < D$ and $l<B^n$.
Consider  Vronskian
\begin{eqnarray}\label{Vronskian}
W(\Phi_1,\ldots,\Phi_l)=\left|\begin{array}{cccc}
\Phi_1(x) & \Phi_2(x) & \ldots & \Phi_l(x) \\
\Phi_1'(x) & \Phi_2'(x) & \ldots & \Phi_l'(x) \\
  \vdots & \vdots & \ddots & \vdots \\
 \Phi_1^{(l-1)}(x) & \Phi_2^{(l)}(x) & \ldots & \Phi_l^{(l-1)}(x)
\end{array}\right|.
\end{eqnarray}
That is a polynomial of $x$ (let us call it $P(x)$) having the degree at most
$$
    \deg P(x)
        \leqslant
            \sum_{i=1}^l \sum_{j=0}^{n-1} t b_{j,i} + l(D-1) -\frac{1}{2}l(l-1) \,.
$$
It is easy to see that $P(x)$ divided by polynomials
$$
    \Psi_0(x)=x^{\sum_{i=1}^l t b_{0,i}-\frac{1}{2}l(l-1)}
$$
and polynomials
$$
\Psi_k(x)=(x-\alpha_k)^{\left(t\sum_{i=1}^l b_{k,i}\right)-\frac{1}{2}l(l-1)},\quad k=1,\ldots,n-1,
$$
which are called $\Psi_0(x),\ldots,\Psi_{n-1}(x)$.
Thus $P(x)$ divided by
$$
\Psi(x)=\prod_{k=0}^{n-1}\Psi_k(x).
$$
At the same time
$$
\deg \Psi(x)=\sum_{i=1}^l tb_{0,i}+t\sum_{k=1}^{n-1}\sum_{i=1}^l b_{k,i}-\frac{1}{2}nl(l-1)
    =
        \deg P(x)-\frac{1}{2}(n-1)l(l-1) - l(D-1) \,.
$$

It is remain to note that if
$P(x)$ divided by $(x-\alpha_n)^{C}$ then $P(x)/\Psi(x)$ divided by the same monomial and
$$
    \deg (P(x)/\Psi(x))\leqslant \frac{1}{2}(n-1)l(l-1) + l(D-1) \,.
$$
Hence either $C\leqslant \frac{1}{2}(n-1)l(l-1)$ or $P(x)\equiv 0$
but in the case the polynomials $\Phi_1(x),\ldots,\Phi_l(x)$ are linearly dependent
(see Lemma \ref{lem-linzav_new} below)
and we reduce the original problem to the question with the smaller number of brackets.

Now return to our suggestion that the sum $\Phi(x)$ from  (\ref{lin-komb2}) divided by $(x-\alpha_n)^t$.
In the case Vronskian $P(x)=W(\Phi_1,\ldots,\Phi_l)$  divided by $(x-\alpha_n)^{t-(l-1)}$
because of the polynomials $\Phi(x),\ldots,\Phi^{(l-1)}(x)$ are divided by $(x-\alpha_n)^{t-(l-1)}$.
Thus
$$
    t \leqslant(l-1)+\frac{1}{2}(n-1)l(l-1) + l (D-1) < \frac{1}{2}(n-1)l(l-1) + lD \,.
$$
On the other hand the total number $l$ of the polynomials $l$ in (\ref{lin-komb2}) is bounded by $l< B^{n}$.
Hence
$$
    t
        <
            \frac{1}{2}(n-1) B^{2n} + D B^n
$$
with contradiction.
This completes the proof.
$\Box$

We give two lemmas on linear independence.
Lemma \ref{lem-linzav} is
a simple
general statement
and Lemma \ref{lem-linzav_new} allows us to have better
dependence between parameters $p$, $t$, $n$, and $B$.

\Lemma
\label{lem-linzav}
{\it
Let Vroskian $P(x)=W(\Phi_1(x),\ldots,\Phi_l(x))$ of degree less than $p$ equals zero in $\mathbb{Z}_p [x]$.
Then there is a nontrivial linear combination
of the polynomials $\Phi_1(x),\ldots,\Phi_l(x)$ with coefficients from $\mathbb{Z}_p$ such that
$$
\mu_1\Phi_1(x)+\ldots+\mu_l\Phi_l(x)\equiv 0,\qquad \mu_1,\ldots,\mu_l\in\mathbb{Z}_p.
$$
}
$\Box$

\Lemma
\label{lem-linzav_new}
{\it
Suppose that the notation of Proposition \ref{pred:polyn} holds.
Let Vroskian $P(x)=W(\Phi_1(x),\ldots,\Phi_l(x))$ equals zero in $\mathbb{Z}_p [x]$
$$
P(x)\equiv 0 \,.
$$
Then there is a nontrivial linear combination of the polynomials
$\Phi_1(x),\ldots,\Phi_l(x)$ with coefficients
$\mu_i\in\mathbb{Z}_p$, $i\in [l]$ such that
$$
\mu_1\Phi_1(x)+\ldots+\mu_l\Phi_l(x)\equiv 0 \,,
$$
provided by $p \ge (2nB+2)t$.
}
\\
\Proof
Since $P(x) \equiv 0$ it follows that
there is a nontrivial zero combination of its rows, i.e.
\begin{eqnarray}\label{lin-ur-syst}
\lambda_1\Phi_k(x)+\lambda_2\Phi_k'(x)+\ldots+\lambda_{l}\Phi_k^{(l-1)}(x)\equiv 0,\qquad k=1,\ldots,l \,,
\end{eqnarray}
where the coefficients $\lambda_i=\lambda_i(x)$ depend on $x$, in general, and does not equal zero simultaneously.
We prove that the coefficients $\lambda_i$ can be chosen
do not depend of $x$
and does not equal zero simultaneously.
Linear combination (\ref{lin-ur-syst}) can be considered as a formal linear differential equation
of the order at most $l-1$:
\begin{eqnarray}\label{lin-ur-vr}
\lambda_1 u(x)+\lambda_2 u'(x)+\ldots+\lambda_{l}u^{(l-1)}(x)=0.
\end{eqnarray}
Polynomials $u(x)$, satisfying the last equation form a linear space.
It is easy to see that any solution of (\ref{lin-ur-vr}) having $l-1$ derivatives at some point $x_0$
equal zero is equal to zero identically.
Indeed, putting, say, $x_0=0$ in (\ref{lin-ur-vr}), we get a linear relation
between $u^{(l-1)}(0)$ and $u^{(l-1)}(0),\ldots,u(0)$.
Taking the formal derivation of (\ref{lin-ur-vr}),
we obtain similar relations for $u^{(l)}(0)$ and so on.
Thus all derivations of $u$ are zero because they can be expressed as linear combinations of $u^{(l-1)}(0),\ldots,u(0)$.
We will prove below that the degrees  of the functions $\lambda_i (x)$
as well as linear combination (\ref{lin-ur-vr})  is less than $p$.
Thus we can take the formal derivations of all these functions and apply the previous arguments.

Now consider a linear combination of columns of the Vronskian at the point $x=0$.
By assumption we have for some $\mu_1,\ldots,\mu_l$ that
\begin{eqnarray*}
\mu_1 \Phi_1^{(k)}(0)+\mu_2\Phi_2^{(k)}(0)+\ldots+\mu_l\Phi_l^{(k)}(0)=0,\qquad k=0,1,\ldots,l-1.
\end{eqnarray*}
Consider the solution
\begin{eqnarray*}
u(x)=\mu_1 \Phi_1(x)+\mu_2\Phi_2(x)+\ldots+\mu_l\Phi_l(x)
\end{eqnarray*}
of equation (\ref{lin-ur-vr}).
Then $u(0), \dots, u^{(l-1)} (0)$  equal zero.
By the previous arguments $u(x) \equiv 0$.
Thus we have found a zero linear combination of the polynomials $\Phi_1(x),\ldots,\Phi_l(x)$
with
coefficients $\mu_1,\ldots,\mu_l \in \Z_p$
and we are done.

It is remain to show that the left hand side of equation (\ref{lin-ur-vr}) is a polynomial of degree less than $p$.

\Lemma
{\it
The degree of the polynomial
$$
\lambda_1 u(x)+\lambda_2 u'(x)+\ldots+\lambda_{l}u^{(l-1)}(x)
$$
less than
$(2nB+2)t$.
}
\\
\Proof The coefficients $\lambda_1,\ldots,\lambda_l$
are solutions of homogeneous system of linear equations (\ref{lin-ur-syst}).
Clearly, system (\ref{lin-ur-syst}) has a nonzero solution for all $x$.
We will use Cramer's rule.
Suppose that there are $l_1$ linear independent equations among $l$ equations of the system.
Without loss of generality one can suppose that these are the first $l_1$ equations.
Further there exist $l_1$ columns of the matrix of system (\ref{lin-ur-syst}) such that the matrix
formed by the elements of the first $l_1$ rows and these $l_1$ columns is non--degenerate for some $x$.
By $i_{1},\ldots,i_{l_1}$ denote the indexes of the columns
and let  $j_1,\ldots,j_{l-l_1}$ be the indexes of another columns.
Let us solve system (\ref{lin-ur-syst}).
We have
$$
\lambda_{i_1}\Phi_k^{(i_1-1)}(x)+\ldots+\lambda_{i_{l_1}}\Phi_k^{(i_{l_1}-1)}(x)=-\sum_{s=1}^{l-l_1}\lambda_{j_s}\Phi_{k}^{(j_s-1)}(x)=\hat{\Phi}_k(x),\qquad
k=1,\ldots,l_1 \,,
$$
where
$$
\hat{\Phi}_k(x)=-\sum_{s=1}^{l-l_1}\lambda_{j_s}\Phi_{k}^{(j_s-1)}(x),\qquad k=1,\ldots,l_1.
$$
The solutions of the system form a linear space of the dimension $l-l_1$.
Put
$\lambda_{j_1},\ldots,\lambda_{j_{l-l_1}}$ equal
$$
\lambda_{j_s}=\frac{x^{D+tB}\prod_{j=1}^{n-1}(x-\alpha_j)^{tB}}{\hat{\Psi}(x)}\left|\begin{array}{ccc}
  \Phi_1^{(i_1-1)}(x) & \ldots & \Phi_{l_1}^{(i_1-1)}(x) \\
  \ldots & \ldots & \ldots \\
  \Phi_1^{(i_{l_1}-1)}(x)& \ldots & \Phi_{l_1}^{(i_{l_1}-1)}(x)
\end{array}\right|,\qquad s=1,\ldots,l-l_1 \,,
$$
where
$$
\hat{\Psi}(x)=x^{\sum_{q=1}^{l_1}(c_q+tb_{0,q})-\sum_{q=1}^{l_1}(i_{l_1}-i_q)}                        
\prod_{j=1}^{n-1}(x-\alpha_j)^{(t\sum_{q=1}^{l_1}b_{j,q})-\sum_{q=1}^{l_1}(i_{l_1}-i_q)} \,,
$$
where $c_q = \deg H_q < D$.
Then by Cramer's rule for $i=1,\ldots,l-1$, we obtain
$$
\lambda_{i_s}=\frac{x^{D+tB}\prod_{j=1}^{n-1}(x-\alpha_j)^{tB}}{\Psi(x)}\left|\begin{array}{ccccccc}
  \Phi_1^{(i_1-1)}(x) & \ldots & \Phi_{l_1}^{(i_1-1)}(x) \\
    \ldots & \ldots & \ldots \\
  \Phi_1^{(i_{s-1}-1)}(x) & \ldots & \Phi_{l_1}^{(i_{s-1}-1)}(x) \\
  \hat{\Phi}^*_1(x) & \ldots & \hat{\Phi}^*_{l_1}(x) \\
   \Phi_1^{(i_{s+1}-1)}(x) & \ldots & \Phi_{l_1}^{(i_{s+1}-1)}(x) \\
       \ldots & \ldots & \ldots \\
  \Phi_1^{(i_{l_1}-1)}(x) & \ldots & \Phi_{l_1}^{(i_{l_1}-1)}(x)
\end{array}\right|,\qquad s=1,\ldots,l_1 \,,
$$
where
$$
    \hat{\Phi}^*_k(x)=-\sum_{s=1}^{l-l_1} \Phi_{k}^{(j_s-1)}(x),\qquad k=1,\ldots,l_1.
$$
It is easy to see that all
$\lambda_1,\ldots,\lambda_l$  are polynomials.
Let us find an upper bound for the degrees of such polynomials
$$
\deg \lambda_i(x)\leqslant \frac{1}{2}l(l-1)(n-1)+(l+1)D+nBt<(nB+1)t,\qquad i=1,\ldots,l;
$$
The degree of each $\Phi_k(x)$
does not exceed
$$
\deg \Phi_k(x)<nBt+D \,,
$$
hence
$$
\deg (\lambda_1\Phi_k(x)+\lambda_2\Phi_k'(x)+\ldots+\lambda_{l}\Phi_k^{(l-1)}(x)) < 2ntB + t + D - 1 < (2nB+2)t \,,
$$
as required.
$\Box$

\Note Proposition \ref{pred:polyn} can be proven using Fuchs equation
for Levelt's basis
(see the formulation in \cite{Bo}).
Nevertheless, we prefer to use a more simple approach calculating the degree of Vronskian of the system
of the polynomials from (\ref{polyn}).

Similarly, we obtain the following proposition.

\Pred
{\it
Let $n$, $t$, $B$, $D$, $D<t$ be positive integers, and  $p$ be a prime number.
Let also $T$ be a set, $T \subseteq \Z_p^*$, $T > n-1$.
Finally, suppose that
\begin{equation}\label{f:polyn_parameters_min}
    t \ge D B^n + \frac{\frac{n}{2} D^2 B^{2n}}{|T|-n+1} \,, \quad
    p \ge \min\{ tnDB^{n+1}, B^2 t^2 n^2 \} \,.
\end{equation}
Then there is a tuple $\a_1, \dots, \a_n \in T$ such that the polynomials of the form
\begin{eqnarray}\label{polyn_min}
    x^{a_i} x^{tb_{0,i}}(x-\alpha_1)^{tb_{1,i}}\ldots (x-\alpha_n)^{tb_{n,i}}
\end{eqnarray}
where $a_i < D$, $b_{0,i},b_{1,i},\ldots,b_{n,i}<B$ are linearly independent over $\Z_p$.
}
\label{pred:polyn_min}
\\
\Proof
One can suppose that for some $\a_1,\dots,\alpha_{n-1} \in T$ the correspondent polynomials from (\ref{polyn_min})
are linearly independent over $\Z_p$, otherwise
we have a problem with smaller number of brackets.
Thus, fix $\a_1,\dots,\alpha_{n-1} \in T$
and let $\a_n$ belongs to the nonempty set $T\setminus \{ \a_1,\dots, \a_{n-1} \}$.
After that apply the arguments as in Proposition \ref{pred:polyn}.
Suppose that there is a nontrivial linear combination of the polynomials from (\ref{polyn_min})
which equals zero identically
\begin{eqnarray}\label{lin-komb'}
\sum_{i=1}^m C_i x^{a_i} x^{tb_{0,i}}(x-\alpha_1)^{tb_{1,i}}\ldots (x-\alpha_n)^{tb_{n,i}}\equiv 0.
\end{eqnarray}
Divide   (\ref{lin-komb'}) by  $(x-\alpha_n)^{ts}$, where $s=\min_{i}b_{n,i}$.
Consider the terms from (\ref{lin-komb'})
with minimal $b_{n,i}$, i.e. equal $s$.
One can suppose that these are the first $l$ terms.
Then the polynomial
\begin{eqnarray}\label{lin-komb2'}
\Phi(x)=\sum_{i=1}^l C_i x^{a_i} x^{tb_{0,i}}(x-\alpha_1)^{tb_{1,i}}\ldots (x-\alpha_{n-1})^{tb_{n-1,i}}
\end{eqnarray}
 divided by $(x-\alpha_n)^t$.
Denote the sum of polynomials from (\ref{lin-komb2'}) by
\begin{eqnarray*}
\Phi_i(x)=x^{a_i} x^{tb_{0,i}}(x-\alpha_1)^{tb_{1,i}}\ldots (x-\alpha_{n-1})^{tb_{n-1,i}},\qquad i=1,\ldots,l.
\end{eqnarray*}
\bigskip
Consider Vronskian
\begin{eqnarray}\label{Vronskian'}
W(\Phi_1,\ldots,\Phi_l)=\left|\begin{array}{cccc}
\Phi_1(x) & \Phi_2(x) & \ldots & \Phi_l(x) \\
\Phi_1'(x) & \Phi_2'(x) & \ldots & \Phi_l'(x) \\
  \vdots & \vdots & \ddots & \vdots \\
 \Phi_1^{(l-1)}(x) & \Phi_2^{(l)}(x) & \ldots & \Phi_l^{(l-1)}(x)
\end{array}\right|.
\end{eqnarray}
That is a polynomial of $x$ (let us call it $P(x)$) having the degree at most
$$
\deg P(x)\leqslant\sum_{i=1}^l\left(a_i+t\sum_{j=0}^{n-1}b_{j,i}\right)-\frac{1}{2}l(l-1).
$$
It is easy to see that $P(x)$ divided by polynomials
$$
\Psi_0(x)=x^{\sum_{i=1}^l (a_i+t b_{0,i})-\frac{1}{2}l(l-1)}
$$
and polynomials
$$
\Psi_k(x)=(x-\alpha_k)^{\left(t\sum_{i=1}^l b_{k,i}\right)-\frac{1}{2}l(l-1)},\quad k=1,\ldots,n-1,
$$
which are called $\Psi_0(x),\ldots,\Psi_{n-1}(x)$.
Thus $P(x)$ divided by
$$
\Psi(x)=\prod_{k=0}^{n-1}\Psi_k(x).
$$
At the same time
$$
    \deg \Psi(x)=\sum_{i=1}^l(a_i+tb_{0,i})+t\sum_{k=1}^{n-1}\sum_{i=1}^l b_{k,i}-\frac{1}{2}nl(l-1)=\deg P(x)-\frac{1}{2}(n-1)l(l-1).
$$

It is remain to note that if $P(x)$ divided by  $(x-\alpha_n)^{C}$ then
 $P(x)/\Psi(x)$  divided by the same monomial and
$$
\deg (P(x)/\Psi(x))\leqslant \frac{1}{2}(n-1)l(l-1) \,.
$$
Hence $C\leqslant \frac{1}{2}(n-1)l(l-1)$.
Note that the polynomial $P(x)$ does not equal zero identically because in the case
the polynomials $\Phi_1(x),\ldots,\Phi_l(x)$ are linearly dependent
and we obtain a contradiction
(see Lemma \ref{lem-linzav} or Lemma \ref{lem-linzav_new}).

Now return to our suggestion that the sum $\Phi(x)$ from (\ref{lin-komb2'})
divided by $(x-\alpha_n)^t$.
In the case Vronskian $P(x)=W(\Phi_1,\ldots,\Phi_l)$ divided by $(x-\alpha_n)^{t-(l-1)}$
because of the polynomials $\Phi(x),\ldots,\Phi^{(l-1)}(x)$ divided by $(x-\alpha_n)^{t-(l-1)}$.
Thus for
{\it every} $\a_n$ Vronskian  $P(x)=W(\Phi_1,\ldots,\Phi_l)$ divided by $(x-\alpha_n)^{t-(l-1)}$.
Whence
$$
    (|T|-n+1) (t-(l-1)) \leqslant(l-1)+\frac{1}{2}(n-1)l(l-1) + l (D-1) < \frac{1}{2}(n-1)l(l-1) + lD \,.
$$
On the other hand the total number of polynomials $l$ in (\ref{lin-komb2'})
equals
$DB^{n}$.
Hence
$$
    t
        <
            D B^n + \frac{\frac{n}{2} D^2 B^{2n}}{|T|-n+1}
$$
with contradiction.

Note also
$$
    \deg P(x) \le \min\{ tn D B^{n+1}, B^2 t^2 n^2 \} < p \,.
$$
This completes the proof. $\Box$

\refstepcounter{section}
\label{sec:proof}

{\bf \arabic{section}. The proof of the main result.}

Let $R\subseteq \Z_p^*$ be a  multiplicative subgroup, and $t=|R|$.
Let $k\ge 1$ be a positive integer, and $\mu_1,\dots,\mu_{k}$
be fixed
different
nonzero elements.
Let also  $\xi_0,\xi_1,\dots,\xi_k$ be some nonzero residuals,
and $A_{\v{\xi},\la}$, $\v{\xi} = (\xi_0,\xi_1,\dots,\xi_k)$, $\la \in \Z_p^*$
be  arbitrary subsets of the set
$\xi_0 R \bigcap (\xi_1 R +\la \cdot \mu_1) \bigcap \dots \bigcap (\xi_k R +\la \cdot \mu_k)$.
Finally, suppose that
we have a family of sets  $A_{\v{\xi}_1,\la_1}, \dots, A_{\v{\xi}_s,\la_s}$,
where
the sets $A_{\v{\xi}_l,\la_l}$
can have the same $\la_l$.

Applying Stepanov's method, we prove one of the main lemmas of the section.
We use arguments from \cite{K_Tula} (see also \cite{Heath_B-K}).

\Lemma
{\it
    Let $R\subseteq \Z_p^*$ be a multiplicative subgroup, and $t=|R|$.
    Let $k \ge 2$, $s$, $B$ be arbitrary positive integers such that
    \begin{equation}\label{cond:B,s,t}
        k B^{2k} < t \,,
                \quad ts < B^{2k+1} \,,
    \end{equation}
    and     \begin{equation}\label{cond:B,s,t_1}
            p \ge (2kB+2) t\,.
    \end{equation}
    Let also $A_{\v{\xi}_1,\la_1}, \dots, A_{\v{\xi}_s,\la_s}$
    be some sets of the family above.
    Then
    \begin{equation}\label{f:E_est}
        \sum_{l=1}^s |A_{\v{\xi}_l,\la_l}| \le \frac{(k+1) t B}{[t/2B^k]} \,.
    \end{equation}
}
\label{l:main_lemma_par}
\Proof
Let $D=[t/(2B^k)]$.
Since $2 B^{k} \le k B^{2k} < t$ it follows that  $D\ge 1$.
Let also $\mathcal{E}$ be the union of all sets $A_{\v{\xi}_l,\la_l}$.
One can assume that the sets  $A_{\v{\xi}_l,\la_l}$ are disjoint, and $\la=1$.
We need to estimate the size of the set $\mathcal{E}$.
Let $\Phi(X,Y,Z_1,\dots,Z_k) \in \Z_p [X,Y,Z_1,\dots,Z_k]$ be an arbitrary polynomial such that
$$
    \deg_X \Phi < D \,, \quad \deg_Y \Phi < B  \,, \quad \deg_{Z_j} \Phi < B \,, \quad  j \in [k] \,.
$$
We have
\begin{equation}\label{f:Phi_pol}
    \Phi(X,Y,Z_1,\dots,Z_k) = \sum_{a,b,\v{c}} \la_{a,b,\v{c}} X^a Y^b Z^{\v{c}} \,,
\end{equation}
where $\v{c} = (c_1,\dots,c_k) \in \Z^k_p$ and $Z^{\v{c}} = Z^{c_1}_1 \dots Z_k^{c_k}$.
Besides
\begin{equation}\label{f:Psi_pol}
    \Psi (X) = \Phi (X,X^t, (X-\mu_1)^t, \dots, (X-\mu_k)^t) \,.
\end{equation}
Clearly
$$
    \deg \Psi \le D-1 + (k+1) t (B-1) \,.
$$
If we will find the coefficients $\la_{a,b,\v{c}}$ such that, firstly,
the polynomial $\Psi$ is nonzero, and, secondly,
$\Psi$ has the root of order at least $D$ at any point of the set $\mathcal{E}$ then
$$
    |\mathcal{E}|
        \le
            (D-1 + (k+1) t (B-1) ) / D
                <
                    \frac{(k+1) t B}{[t/2B^k]}
$$
and lemma will be proved.
Thus, we should check that
$$
    \l( \frac{d}{d X} \r)^n \Psi (X) \Big|_{X=x} = 0 \,, \quad \forall n < D \,, \quad \forall x\in \mathcal{E} \,.
$$
For any $x\in \mathcal{E}$, we have $x\neq 0$ and $x\neq \mu_j$, $j\in [k]$.
Hence the last condition is equivalent
\begin{equation}\label{f:want_Psi}
   [X (X-\mu_1) \dots (X-\mu_k)]^n \l( \frac{d}{d X} \r)^n \Psi (X) \Big|_{X=x} = 0 \,, \quad \forall n < D \,, \quad \forall x\in \mathcal{E} \,.
\end{equation}
It is easy to see that for all $m,q$, $q\ge m$, and any $\mu$ the following holds
$$
    (X-\mu)^m \l( \frac{d}{d X} \r)^m (X-\mu)^q = \frac{q!}{(q-m)!} (X-\mu)^q \,.
$$
If $m>q$ then  the left hand side equals zero.
So, there are well--defined polynomials  $P_{n,a,b,\v{c}} (X)$ such that
$$
    [X (X-\mu_1) \dots (X-\mu_k)]^n \l( \frac{d}{d X} \r)^n X^a X^{tb} (X-\mu_1)^{t c_1} \dots (X-\mu_k)^{t c_k}
        =
$$
$$
        =
           P_{n,a,b,\v{c}} (X) X^{tb} (X-\mu_1)^{t c_1} \dots (X-\mu_k)^{t c_k} \,.
$$
Here $a,b,c_1,\dots,c_k$ are nonnegative integers.
For some $a,b,\v{c}$ polynomial $P_{n,a,b,\v{c}}$ can be identically zero.
Clearly,  $\deg P_{n,a,b,\v{c}} \le a+n$.
By the definition of the sets
$\mathcal{E}$ and $A_{\v{\xi}_l,\v{\mu}_l}$, we have
$$
    [X (X-\mu_1) \dots (X-\mu_k)]^n \l( \frac{d}{d X} \r)^n X^a X^{tb} (X-\mu_1)^{t c_1} \dots (X-\mu_k)^{t c_k}  \Big|_{X=x}
        =
$$
$$
        =
            y^b_0 (l) y^{c_1}_1 (l) \dots y^{c_k}_1 (l) P_{n,a,b,\v{c}} (X) \,, \quad x\in A_{\v{\xi}_l,\v{\mu}_l} \,,
$$
where residuals $y^b_0 (l), y^{c_1}_1 (l), \dots, y^{c_k}_1 (l)$ does not depend on
the choice of the element $x\in A_{\v{\xi}_l,\v{\mu}_l}$.
By (\ref{f:Phi_pol}),  (\ref{f:Psi_pol})
$$
    [X (X-\mu_1) \dots (X-\mu_k)]^n \l( \frac{d}{d X} \r)^n \Psi (X) \Big|_{X=x}
        =
$$
$$
        =
            \sum_{a,b,\v{c}} \la_{a,b,\v{c}}\, \cdot y^b_0 (l) y^{c_1}_1 (l) \dots y^{c_k}_1 (l)
                P_{n,a,b,\v{c}} (x) := P_{n,\,l} (x) \,, \quad x\in A_{\v{\xi}_l,\v{\mu}_l} \,.
$$
Coefficients of the  polynomials $P_{n,\,l} (X)$ are linear forms of $\la_{a,b,\v{c}}$.
Choose $\la_{a,b,\v{c}}$ such that polynomials $P_{n,\,l} (X)$ are identically zero
for an arbitrary $n<D$ and any $l\in [s]$.
Then equality (\ref{f:want_Psi}) holds for all $x\in \mathcal{E}$.
We have (\ref{cond:B,s,t}).
Since $\deg P_{n,\,l} < 2D$ it follows that
\begin{equation}\label{f:main_connection_between_parameters}
    2sD^2 \le 2s Dt/2B^k < D B^{k+1} \,,
\end{equation}
and
(\ref{f:main_connection_between_parameters})
guarantee that there is a nonzero tuple of coefficients $\la_{a,b,\v{c}}$
such that
$P_{n,\,l} (X) \equiv 0$, $n<D$, $l\in [s]$.

We must check that the obtained
polynomial $\Psi (X)$ is nonzero.
We have  $D=[t/(2B^k)]$, and $k B^{2k} < t$.
Hence $t \ge \frac{1}{2}(k-1) B^{2k} + D B^k$.
Besides inequality  (\ref{cond:B,s,t_1}) holds.
Using Proposition \ref{pred:polyn} with $n=k$, we obtain
that the polynomial $\Psi (X)$ is nonzero identically.
This concludes the proof of the lemma. $\Box$

{\bf Proof of Theorem \ref{t:main_many_shifts}.}
Let $|R|=t$, $s=|Q|/t$.
Let 
$B$ be the least integer such that $B^{2k+1} > ts$.
Then $B \le (ts)^{1/(2k+1)} + 1$.
Using bound
$$
    |Q|
        < ( (t/k)^{1/2k} - 1)^{2k+1} \,,
$$
we get
$$
    k B^{2k} \le k ((ts)^{1/(2k+1)} + 1)^{2k} < t
$$
and  condition (\ref{cond:B,s,t}) of Lemma  \ref{l:main_lemma_par} is satisfied.
Since  $t> k 2^{2k+4}$ and
$$
    |Q| < ( (t/k)^{1/2k} - 1)^{2k+1} < (t/k)^{(2k+1)/2k} \,,
$$
it follows that  $t/2B^{k} \ge 2$.
Finally, inequality  (\ref{cond:B,s,t_1}) of the same Lemma is a consequence
of $p\ge 4tk ( |Q|^{\frac{1}{2k+1}} + 1)$.
Applying the lemma and using the bounds $t/2B^{k} \ge 2$, $B \le |Q|^{\frac{1}{2k+1}} + 1$, we obtain
$$
    \sum_{\la \in Q} |R \bigcap (R +\la \cdot \mu_1) \bigcap \dots \bigcap (R +\la \cdot \mu_k)|
        \le
            t \frac{(k+1) t B}{[t/2B^k]}
                \le
                    4 (k+1) B^{k+1} t
                        \le
                            4 (k+1) (|Q|^{\frac{1}{2k+1}} + 1)^{k+1} t \,.
$$
This completes the proof.
$\Box$

{\bf Proof of Corollary  \ref{cor:main_many_shifts}.}
It is sufficiently to check that for all $|R| \ge 32 k 2^{20k \log (k+1)} > k 2^{2k+4}$ the following holds
$$
    |R| < ( (|R|/k)^{1/2k} - 1)^{2k+1} \,.
$$
It is easy to see that the assumed bounds for the cardinality of $R$ imply the last inequality.
$\Box$

Using Proposition  \ref{pred:polyn_min} instead of Proposition  \ref{pred:polyn},
we obtain the following statement.

\St
{\it
    Let $k \ge 2$ be a positive integer,
    and $R\subseteq \Z_p^*$ be a multiplicative subgroup.
    Let also $T\subseteq \Z_p^*$ be any set, $2k \le |T| \le |R|k /2$,
    and let $s$, $B$ be arbitrary natural numbers such that
    \begin{equation}\label{cond:B,s,t_min}
         2k B^{2k} \le |R| |T|,
                \quad 2s \l( \frac{|R||T|}{2k} \r)^{1/2} <  B^{2k+1} \,,
    \end{equation}
    and
    \begin{equation}\label{cond:B,s,t_1_min}
            p \ge (2kB+2)t\,.
    \end{equation}
    Then there are different elements $\mu_j \in T$, $j\in [k]$
    such that
    for all sets $A_{\v{\xi}_1,\la_1}, \dots, A_{\v{\xi}_s,\la_s}$
    the following holds
    \begin{equation}\label{f:E_est_min}
        \sum_{l=1}^s |A_{\v{\xi}_l,\la_l}| \le \frac{(k+1) |R| B}{[( |R||T|/(2k B^{2k}) )^{1/2}]} \,.
    \end{equation}
}
\label{st:main_lemma_par_min}
\Proof
Let $t=|R|$, and
$D=[( t|T|/(2k B^{2k}) )^{1/2}]$.
Since
$2k B^{2k} \le t |T|$ it follows that  $D\ge 1$.
Besides $D<t$ because of $|T| \le tk/2$.
Let also $\mathcal{E}$ be the union of all sets $A_{\v{\xi}_l,\la_l}$.
Using the arguments as in Lemma  \ref{l:main_lemma_par}, we construct a polynomial $\Psi$,
having a root of order at least $D$ at any point of the set $\mathcal{E}$.
If the polynomial $\Psi$ is nonzero then we have the following bound for the cardinality of the set $\mathcal{E}$
$$
    |\mathcal{E}|
        \le
            (D-1 + (k+1) t (B-1) ) / D
                <
                    \frac{(k+1) t B}{[( t|T|/(2k B^{2k}) )^{1/2}]} \,.
$$
Besides an analog of inequality (\ref{f:main_connection_between_parameters})
is
\begin{equation}\label{f:main_connection_between_parameters'}
    2sD^2 \le 2s D ( t|T|/(2k B^{2k}) )^{1/2} < D B^{k+1} \,,
\end{equation}
where the second inequality from (\ref{cond:B,s,t_min})
was used.
By (\ref{cond:B,s,t_min}) and  $|T| \le tk/2$, we find that
$$
    t \ge D B^k + \frac{\frac{k}{2} D^2 B^{2k}}{|T|-k+1} \,.
$$
Using condition (\ref{cond:B,s,t_1_min}) and applying Proposition  \ref{pred:polyn_min} with $n=k$,
we obtain that
for some different $\mu_j \in T$, $j\in [k]$,
the polynomial
$\Psi (X)$ is nonzero identically.
That concludes the proof.
 $\Box$


\Note
Though sum (\ref{f:E_est_min}) in Statement \ref{st:main_lemma_par_min}
considered for specific tuple of elements $\mu_j$
the dependence between parameters $B,t$ and $T$
(see the first inequality from (\ref{cond:B,s,t_min}))
not so onerousness as bound (\ref{cond:B,s,t}) of Lemma  \ref{l:main_lemma_par}.

\Cor
{\it
    Let $R\subseteq \Z_p^*$ be a multiplicative subgroup, $k\ge 1$ a positive integer.
    Let also $T\subseteq \Z_p^*$ be any set, $2k \le |T| \le |R|k/2$,
    $Q=RQ$ be a  $R$---invariant set,
    $0\notin Q$,
    \begin{equation}\label{f:|Q|_est}
        |Q| < \l( \frac{k |R|}{2|T|} \r)^{1/2} \l( \l( \frac{|R||T|}{8k} \r)^{1/2k} - 1 \r)^{2k+1}
    \end{equation}
    and     
    \begin{equation}\label{f:p_est}
        p \ge \l( \frac{k|R|^3 |T|}{2} \r)^{1/2} \l( \l( |Q| \l( \frac{2|T|}{k|R|} \r)^{1/2} \r)^{1/(2k+1)} + 1 \r) \,.
    \end{equation}
    Then
    $$
        \min_{\mu_1,\dots,\mu_k \in T,\, \mu_i \neq \mu_j}\, C_{k+1} (Q,R,\dots,R) (\mu_1,\dots,\mu_k)
            \le
    $$
    $$
            \le
                    (32 k^3)^{1/2} \l( \frac{|R|}{|T|} \r)^{1/2}
                \l( \l( |Q| \l( \frac{2|T|}{k|R|} \r)^{1/2} \r)^{1/(2k+1)} + 1 \r)^{k+1} \,.
    $$
}
\label{cor:C_k(Q,R)}
\Proof
Let $t=|R|$, $s=|Q|/t$.
Let $B$ be the least integer such such that $B^{2k+1} > 2s \l( \frac{t|T|}{2k} \r)^{1/2}$.
Then $B\le \l( 2s \l( \frac{t|T|}{2k} \r)^{1/2} \r)^{1/(2k+1)} + 1$.
Since  $|Q| < \l( \frac{k t}{2|T|} \r)^{1/2} \l( \l( \frac{|R||T|}{8k} \r)^{1/2k} - 1 \r)^{2k+1}$
it follows that
\begin{equation}\label{tmp:27_12_10_1}
    2k B^{2k} \le 8k B^{2k} \le 8k \l( \l( 2s \l( \frac{t|T|}{2k} \r)^{1/2} \r)^{1/(2k+1)} + 1 \r)^{2k} < t |T| \,.
\end{equation}
Thus all conditions (\ref{cond:B,s,t_min})
of Statement \ref{st:main_lemma_par_min}
are satisfied.
Inequality (\ref{cond:B,s,t_1_min}) of the lemma is a consequence of bound (\ref{f:p_est}).
Using 
Statement \ref{st:main_lemma_par_min}
and  (\ref{tmp:27_12_10_1}), we obtain
$$
    \min_{\mu_1,\dots,\mu_k \in T}\, C_{k+1} (Q,R,\dots,R) (\mu_1,\dots,\mu_k)
        \le
            \frac{(k+1) t B}{[( t|T|/(2k B^{2k}) )^{1/2}]}
                \le
                    (8 k)^{1/2} (k+1) \l( \frac{t}{|T|} \r)^{1/2} B^{k+1}
$$
$$
                        \le
                            (32 k^3)^{1/2} \l( \frac{t}{|T|} \r)^{1/2}
                                \l( \l( |Q| \l( \frac{2|T|}{k|R|} \r)^{1/2} \r)^{1/(2k+1)} + 1 \r)^{k+1} \,.
$$
This completes the proof.
$\Box$

\Note One can generalize Corollary \ref{cor:C_k(Q,R)} consider the sum
$\sum_{\la \in Q_1} C_{k+1} (Q,R,\dots,R) (\mu_1,\dots,\mu_k)$, $Q_1 = R Q_1$ as in Theorem \ref{t:main_many_shifts}.
We do not need in the generalization.


\refstepcounter{section}
\label{sec:applications}

{\bf \arabic{section}. On subgroups sumsets.}

First of all we write simple consequences of Lemma \ref{l:main_lemma_par}
and Theorem \ref{t:main_many_shifts}.

\Cor
{\it
    Let $R\subseteq \Z^*_p$ be a multiplicative subgroup, and
    $Q,Q_1,Q_2\subseteq \Z^*_p$ be arbitrary $R$--invariant sets.
    Then \\
    $1)~$ If $|Q| \ll |R|^3$, $|Q| |R|^3 \ll p^3$ then
        \begin{equation}\label{f:improved_Konyagin_old1}
            \sum_{x\in Q} (R\circ R) (x) \ll |R| |Q|^{2/3} \,.
        \end{equation}
    $2)~$ If $|Q| |Q_1| \ll |R|^4$, $|Q| |Q_1| |R|^2 \ll p^3$ then
        \begin{equation}\label{f:improved_Konyagin_old2}
            \sum_{x\in Q} (Q_1\circ R) (x) \ll |R|^{1/3} (|Q||Q_1|)^{2/3} \,.
        \end{equation}
    $3)~$ If $|Q| |Q_1| |Q_2| \ll |R|^{5}$, $|Q| |Q_1| |Q_2| |R| \ll p^3$ then
        \begin{equation}\label{f:improved_Konyagin_old3}
            \sum_{x\in Q} (Q_1 \circ Q_2) (x) \ll |R|^{-1/3} (|Q||Q_1||Q_2|)^{2/3} \,.
        \end{equation}
}
\label{cor:improved_Konyagin_old}

\Note
Clearly,  inequality  (\ref{f:improved_Konyagin_old3}) can be improved
provided that some information
of the set
$\bigcup_{q\in Q} (q^{-1} Q_1 \m q^{-1} Q_2)$
(which is a  multiplicative analog of the set from (\ref{def:otimes})) is known.

Corollary  \ref{cor:improved_Konyagin_old}
implies a statement about additive energy of any $R$---invariant set.
The statement is a tiny generalization of a result from \cite{K_Tula}.
Applying Statement \ref{st:E(Q)_inv} below it is easy to obtain
(using Lemma \ref{l:KB_high_convolutions}, for example)
that
any $R$--invariant set $Q\subseteq \Z^*_p$,
such that $|Q| \ll |R|^{3/2}$, $|Q| |R|^{1/2} \ll p$ has the extension property, namely,
$|Q \pm Q| \gg |Q| |R|^{1/2}$.

\St
{\it
    Let $R\subseteq \Z^*_p$ be a  multiplicative subgroup,
    and $Q \subseteq \Z^*_p$ be an arbitrary $R$--invariant set .
    Let also $|Q| \ll |R|^{3/2}$, and $|Q| |R|^{1/2} \ll p$.
    Then
    \begin{equation}\label{f:E(Q)_inv}
        E(Q) \ll \frac{|Q|^3}{|R|^{1/2}} \quad \mbox{ and } \quad \max_{\xi\neq 0} |\F{Q} (\xi)| \ll |Q|^{7/8} |R|^{-1/4} p^{1/8} \,.
    \end{equation}
}
\label{st:E(Q)_inv}
\Proof
Let $a$ be a parameter.
We have
$$
    E(Q) \le a |Q|^2 + \sum_{x ~:~ (Q\circ Q) (x) \ge a} (Q\circ Q)^2 (x) \,.
$$
Let us arrange values $(Q\circ Q) (x)$, $x\in \Z_p /R$ in decreasing order
and denote its values as $N_1 \ge N_2 \ge \dots$.
Using formula (\ref{f:improved_Konyagin_old3}) of Corollary  \ref{cor:improved_Konyagin_old},
we get $N_j \ll |Q|^{4/3} |R|^{-2/3} j^{-1/3}$.
Hence
$$
    E(Q)
        \ll
            a |Q|^2 + |R| \sum_{j ~:~ j\ll |Q|^4 /(|R|^2 a^3)} j^{-2/3} \cdot \frac{|Q|^{8/3}}{|R|^{4/3}}
                \ll
    a |Q|^2 + \frac{|Q|^4}{|R|a} \,.
$$
Putting $a=|Q|/|R|^{1/2}$, we obtain the required result.
The second inequality in (\ref{f:E(Q)_inv}) is a consequence of the first one, see e.g.
the proof of Corollary  2.5 from \cite{Schoen_Shkr}. $\Box$

We need in a lemma from \cite{Schoen_Shkr}.

\Lemma
{\it
    Let $R\subseteq \Z^*_p$ be a multiplicative subgroup, $|R| \ll p^{2/3}$.
    Then
    $$
        E_3 (R) \ll |R|^3 \log |R| \,.
    $$
}
\label{l:E_3(R)}

Let us obtain a new result on doubling constant of multiplicative subgroups.

\Th
{\it
    Let $R\subseteq \Z^*_p$ be a multiplicative subgroup.
    If $|R| \ll p^{1/2}$ then
    \begin{equation}\label{f:new_R+R}
        |R\pm R| \gg \frac{ |R|^{5/3}}{\log^{1/2} |R|} \,.
    \end{equation}
}
\label{t:subgroups_doubling}
\Proof
Let $S=(R-R)\setminus \{ 0 \}$ (for  $R+R$ we use similar arguments).
Using Lemma \ref{l:KB_high_convolutions} and Corollary \ref{cor:card_KK}, we get
$$
    |R|^6 \le E_3(R) \cdot \sum_{x\in S} (S\circ S) (x) \,.
$$
If $|S| \gg |R|^{5/3}$ then  it is nothing to prove.
In the opposite case, we have $|S|^3 |R| \ll p^3$,
because of the assumption $|R| \ll p^{1/2}$.
Using bound (\ref{f:improved_Konyagin_old3}) of Corollary  \ref{cor:improved_Konyagin_old}
with $Q=Q_1=Q_2=S$,
and Lemma \ref{l:E_3(R)}, we get
$$
    |R|^6 \ll |R|^3 \log |R| \cdot |S|^2 |R|^{-1/3} \,.
$$
Hence $|S| \gg |R|^{5/3} \log^{-1/2} |R|$.
Theorem is proved.
$\Box$

Inequality (\ref{f:new_R+R})
answered on a question of article \cite{Waring_Z_p}.
A weaker bound for subgroups such that $|R| \ll \sqrt{p}$,
better than (\ref{f:HK_3/2}) was obtained by T. Schoen and the second author in \cite{Schoen_Shkr}.
The strongest result on the cardinality of $R\pm R$, where $\sqrt{p} \ll |R| \ll p^{2/3}$,
is contained in \cite{Schoen_Shkr}.
Let us note a consequence of the theorem above.

\Cor
{\it
    Let $R\subseteq \Z^*_p$ be a multiplicative subgroup,
    and $\kappa > \frac{33}{67}$ be a real number.
    Suppose that  $|R| \ge p^{\kappa}$.
    Then for all sufficiently large $p$ the following holds $\Z_p^* \subseteq 6R$.
}
\label{cor:6R_new}


Corollary \ref{cor:6R_new} is a consequence of Theorem \ref{t:subgroups_doubling}
and can be proved exactly as Theorem 4.1 from \cite{Schoen_Shkr},
where the inclusion $\Z_p^* \subseteq 6R$ was obtained under the  assumption $\kappa > \frac{41}{83}$.
Note that a result of  A.A. Glibichuk \cite{Glibichuk_zam} (see also \cite{Rudnev})
implies that
$|4R| > p/2$ (and hence $8R = \Z_p$), provided by $|R| > \sqrt{p}$.

\end{document}